\documentclass[11pt, reqno]{amsart}

\author[P.~Leonetti]{Paolo Leonetti}
\address{Universit\`a degli Studi dell'Insubria, via Monte Generoso 71, 21100 Varese, Italy}
\email{leonetti.paolo@gmail.com}
\urladdr{\url{https://sites.google.com/site/leonettipaolo/}} 



\title{Nonpathological submeasures and $F_\sigma$ ideals}

\usepackage[T1]{fontenc}
\usepackage{amsmath}
\usepackage{amssymb}
\usepackage{amsthm}
\usepackage[left=3.5cm, right=3.5cm]{geometry}
\usepackage{hyperref}
\usepackage{fancyhdr}
\usepackage{enumitem}
\usepackage{comment}
\usepackage{nicefrac}
\usepackage{comment}
\usepackage{bm}
\usepackage{mathrsfs}
\usepackage{graphicx}
\usepackage[utf8]{inputenc}
\usepackage{cancel}
\usepackage{mathtools}

\usepackage{tikz}
\usetikzlibrary{decorations.pathreplacing,matrix,arrows,positioning,automata,shapes,shadows,calc,fadings,decorations,snakes,through,intersections}

\AtBeginDocument{%
   \def\MR#1{}
}

\newtheorem{thm}{Theorem}[section]
\newtheorem{cor}[thm]{Corollary}
\newtheorem{lem}[thm]{Lemma}
\newtheorem{prop}[thm]{Proposition}

\theoremstyle{definition} 
\newtheorem{defi}[thm]{Definition}
\let\olddefi\defi
\renewcommand{\defi}{\olddefi\normalfont}
\newtheorem{question}{Question}
\let\oldquestion\question
\renewcommand{\question}{\oldquestion\normalfont}
\newtheorem{example}[thm]{Example}
\let\oldexample\example
\renewcommand{\example}{\oldexample\normalfont}
\newtheorem{rmk}[thm]{Remark}
\let\oldrmk\rmk
\renewcommand{\rmk}{\oldrmk\normalfont}
\newtheorem{claim}{\textsc{Claim}}

\hypersetup{
    pdftitle={Nonpathological submeasures and $F_\sigma$ ideals},
    pdfauthor={},
    pdfmenubar=false,
    pdffitwindow=true,
    pdfstartview=FitH,
    colorlinks=true,
    linkcolor=blue,
    citecolor=green,
    urlcolor=cyan
}

\providecommand{\MR}[1]{}
\providecommand{\bysame}{\leavevmode\hbox to3em{\hrulefill}\thinspace}
\providecommand{\MR}{\relax\ifhmode\unskip\space\fi MR }

\providecommand{\href}[2]{#2}

\begin{document}

\maketitle
\thispagestyle{empty}

\begin{abstract}
\noindent We study the relationship between nonpathological lower semicontinuous submeasures, $F_\sigma$ ideals on $\omega$, and matrix summability ideals. First, we construct in $\mathsf{ZFC}$ a non-analytic ideal on $\omega$ which is representable as an intersection of matrix summability ideals. Then, we show that an $F_\sigma$ ideal is nonpathological if and only if it can be represented as an intersection of summable ideals if and only if it can be represented as an intersection of matrix summability ideals. We also prove that the random graph ideal is pathological. 
Finally, we compute the exact degree of pathology of a certain submeasure related to a construction of Topsøe.
\end{abstract}

\section{Introduction and Main results}

A map $\varphi: \mathcal{P}(\omega)\to [0,\infty]$ is said to be \emph{submeasure} if it is monotone, subadditive, and $\varphi(\emptyset)=0$. If, in addition, $\varphi$ is finitely additive, that is, $\varphi(A\cup B)=\varphi(A)+\varphi(B)$ for all disjoint $A,B\subseteq X$, then $\varphi$ is a finitely additive measure. (As usual, countably infinite sets are identified with $\omega$ through fixed bijections.) For each submeasure $\varphi$, define
$$
\widehat{\varphi}(B):=\mathrm{sup}\left\{\mu(B): \mu \text{ is finitely additive measure with }\mu\le \varphi\right\},
$$
for each $B\subseteq \omega$. The submeasure $\varphi$ is said to be \emph{nonpathological} if $\varphi=\widehat{\varphi}$, see \cite{FTpathology} and references therein. In addition, the degree of pathology of $\varphi$ is defined by 
$$
P(\varphi):=\mathrm{sup}\left\{\frac{\varphi(B)}{\widehat{\varphi}(B)}: B\subseteq \omega\right\},
$$
with the convention that $0/0:=\infty/\infty:=1$ and $a/0:=\infty/a:=\infty$ for all $a \in (0,\infty)$. 
In addition, a submeasure $\varphi$ is \emph{lower semicontinuous}, shortened as \emph{lscsm}, if $\varphi(A)=\lim_n \varphi(A\cap n)$ for all $A\subseteq \omega$, where each $n \in \omega$ is identified with $\{0,1,\ldots,n-1\}$.

Let $\mathcal{I}$ be an ideal on $\omega$, that is, a hereditary family of subsets of $\omega$ which is stable under finite unions. Unless otherwise stated, we assume that $\omega \notin \mathcal{I}$ and that $[\omega]^{<\omega} \subseteq \mathcal{I}$. An ideal $\mathcal{I}$ on $\omega$ is said to be a $P$-ideal if for every sequence $(S_n)$ in $\mathcal{I}$ there exists $S \in \mathcal{I}$ such that $S_n\setminus S$ is finite for all $n \in \omega$. Identifying $\mathcal{P}(\omega)$ with the Cantor space $\{0,1\}^\omega$, we can speak about the topological complexity of ideals, see e.g. \cite{MR1711328}. 
Recall that an infinite nonnegative matrix $A=\left(a_{n,k}: n,k\in\omega\right)$ is said to be \emph{regular} if it maps convergent real sequences into convergent real sequences, preserving their corresponding limits. By the well-known Silverman--Toepliz theorem, a nonnegative matrix $A$ is regular if and only if 
$$
\lim_{n\to \infty}\sum_{k\in\omega}a_{n,k}=1
\quad \text{ and }\quad 
\lim_{n\to\infty}a_{n,k}=0
\,\,\text{for all }k\in\omega,
$$
see e.g. \cite[Theorem 1.3.3]{MR225045}; cf. also  \cite{MR5065056} and references therein. 
Following \cite{MR1865750}, an ideal $\mathcal{I}$ on $\omega$ is said to be a \emph{matrix summability ideal} if there exists a nonnegative regular matrix $A=(a_{n,k}: n,k \in \omega)$ such that 
$$
\mathcal{I}=\mathcal{I}_A:=\left\{S\subseteq \omega: \lim_{n\to \infty} \sum_{k \in S}a_{n,k}=0\right\}.
$$
Detailed studies on the relationships between ideal convergence and summability can be found in \cite{Filipow18, MR4041540, MR4756411}. In a remarkable result, Jacek Tryba proved in \cite[Corollary 3.12]{MR4756411} that an analytic $P$-ideal $\mathcal{I}$ on $\omega$ can be represented as intersection of matrix summability ideals if and only if it is \emph{nonpathological}, that is, there exists a nonpathological lscsm $\varphi$ such that 
$$
\mathcal{I}=\mathrm{Exh}(\varphi):=\left\{S\subseteq \omega: \lim_{n\to \infty} \varphi(S\setminus n)=0\right\}.
$$
Nonpathological analytic $P$-ideals appeared also in the study of the geometry of Banach spaces, see \cite{MR3436368} and \cite[Section 5.2]{MR4797308}; cf. also \cite{MR2743987}.  
In addition, it has been recently proved in \cite[Corollary 14.4]{FTpathology} that it is consistent in \textsf{ZFC} that there exists 
a non-Borel ideal on $\omega$ which can be written as intersection of matrix summability ideals.  
In the following result, we improve it by showing in \textsf{ZFC} that the latter ideal might be non-analytic. In particular, this answers \cite[Question 3]{Filipow18}. 
\begin{thm}\label{thm:bigcapmatrix} 
There exists a non-analytic ideal $\mathcal{I}$ on $\omega$ which can be represented as intersection of matrix summability ideals. 
\end{thm}

Since every matrix summability ideal is an analytic $P$-ideal (hence, meager), see e.g. \cite{FTpathology}, it follows that the ideal constructed in Theorem \ref{thm:bigcapmatrix} is meager as well.

Next, recall, by a classical result of Mazur \cite{MR1124539}, that an ideal $\mathcal{I}$ on $\omega$ is $F_\sigma$ if and only if there exists a lscsm $\varphi: \mathcal{P}(\omega)\to [0,\infty]$ such that 
$$
\mathcal{I}=\mathrm{Fin}(\varphi):=\{S\subseteq \omega: \varphi(S)<\infty\}
$$
and $\varphi(\{n\})<\infty$ for all $n \in \omega$. 
Thus, following \cite[p. 53]{MR1711328}, we say that an $F_\sigma$ ideal $\mathcal{I}$ on $\omega$ is \emph{nonpathological} if there exists a nonpathological lscsm $\varphi: \mathcal{P}(\omega)\to [0,\infty]$ such that $\mathcal{I}=\mathrm{Fin}(\varphi)$; 
cf. also \cite[Definition 8.3]{FTpathology} and Corollary \ref{cor:nonpathological} below. 
Among nonpathological $F_\sigma$ $P$-ideals, we recall the summable ideals, that is, families of the type $\mathcal{I}_h:=\{S\subseteq \omega: \sum_{n}h(n)<\infty\}$ for some function $h: \omega \to (0,\infty)$ with $\sum_n h(n)=\infty$. 

The next result characterizes the case where an $F_\sigma$ ideal $\mathcal{I}$ is nonpathological. This provides an answer to \cite[Question 11.6]{FTpathology}. 
\begin{thm}\label{thm:representation}
Let $\mathcal{I}$ be an $F_\sigma$ ideal on $\omega$. Then the following are equivalent\textup{:}
\begin{enumerate}[label={\rm (\roman*)}]
\item \label{item:1repres} $\mathcal{I}$ is nonpathological\textup{;}
\item \label{item:2repres} $\mathcal{I}$ can be represented as intersection of summable ideals\textup{;}
\item \label{item:3repres} $\mathcal{I}$ can be represented as intersection of matrix summability ideals\textup{.}
\end{enumerate}
\end{thm}

As a consequence of Theorem \ref{thm:representation}, we show that the definitions of nonpathological ideals in the realms of $F_\sigma$ ideals and analytic $P$-ideals are consistent with each other. 
\begin{cor}\label{cor:nonpathological}
    Let $\mathcal{I}$ be an $F_\sigma$ $P$-ideal on $\omega$. Then $\mathcal{I}$ is a nonpathological $F_\sigma$ ideal if and only if it is a nonpathological analytic $P$-ideal. 
\end{cor}

Now, let us recall that Mazur provided in \cite{MR1124539} the first example of an $F_\sigma$ ideal $\mathcal{M}$ which cannot be extended to a summable one (hence, $\mathcal{M}$ is pathological by Theorem \ref{thm:representation}). 
Recently, Mart{\'i}nez, Meza-Alc{\'a}ntara, and Uzc{\'a}tegui proved in \cite{MR4797308} that in fact $\mathcal{M}$ has a stronger pathological property, namely, every $F_\sigma$ ideal $\mathcal{I}$ containing $\mathcal{M}$ is pathological as well; cf. \cite[Section 5.1]{Carlos26}. In our next result, we are going to show that the \emph{random graph ideal} $\mathcal{R}$ is pathological, which answers \cite[Question 5.19]{Carlos26}. 

To recall the definition of $\mathcal{R}$, pick an independent family $\{W_n: n \in \omega\}$ such that $n \in W_m$ if and only if $m \in W_n$, and consider the random graph $(\omega, E)$ where $\omega$ is the set of vertices and $E:=\{\{n,m\} \in [\omega]^2: n \in W_m\}$ is the set of edges. A set $S\subseteq \omega$ is said to be a clique if $[S]^2 \subseteq E$; also, $S\subseteq \omega$ is said to be an anticlique if $[S]^2\cap E=\emptyset$. Accordingly, $\mathcal{R}$ is the ideal on $\omega$ generated by the family of cliques and anticliques. It is known that $\mathcal{R}$ is a tall $F_\sigma$ ideal, see e.g. \cite[pp. 29-30]{Mezathesis}. 

\begin{thm}\label{thm:randomgraph}
The random graph ideal $\mathcal{R}$ is a pathological $F_\sigma$ ideal. 
\end{thm}

In our last result, we calculate an explicit degree of pathology. To be more precise, recall that Herer and Christensen proved in \cite[Theorem 2]{MR412369} that there exists a nonzero submeasure $\varphi$ such that $\widehat{\varphi}=0$, $\varphi(\omega)=1$, and $\varphi(F)=0$ for all finite $F\subseteq \omega$. Shortly after, Topsøe provided an explicit construction of a submeasure with these properties in \cite[Example 3]{MR419712}. Below, we recall his construction. 

For each integer \(n\geq 1\) and each subset $I\subseteq \Delta_n:=\{1,2,\dots,2^n\}$, define 
$$
A(n,I)
=
\{x\subseteq \Delta_n: |x|=2^{n-1} \text{ and }x \cap I\neq \emptyset\}.
$$
Set $X_n:=A(n,\Delta_n)$. 
For each $n\ge 1$, let $\chi_n: \mathcal{P}(X_n)\to [0,\infty]$ be the submeasure defined by 
$$
\chi_n(B)
=
\min\left\{
\frac{|I|}{2^{n-1}}
:
I \subseteq \Delta_n \text{ and }
B\subseteq A(n,I)
\right\}
$$
for each $B\subseteq X_n$. 
At this point, define also the countably infinite set 
$X:=\bigcup_{n\geq 1}X_n$ 
and let $\chi: \mathcal{P}(X)\to [0,\infty]$ be the submeasure defined by 
\begin{equation}\label{eq:definitionchi}
\chi(B):=\limsup_{n\to \infty}\chi_n(B\cap X_n)
\end{equation}
for each $B\subseteq X$. Then the submeasure $\chi$ has the pathological properties claimed in \cite[Theorem 2]{MR412369}. In fact, it follows that $P(\chi)=\infty$. 

Lastly, let also $\delta: \mathcal{P}(X)\to [0,\infty]$ be the submeasure defined by 
\[
\delta(B)=
\begin{cases}
0 & \text{ if } B=\varnothing,\\
1 & \text{ if }B\neq\varnothing,
\end{cases}
\]
for each $B\subseteq X$. Accordingly, we provide an answer to \cite[Question 7.11]{FTpathology}.
\begin{thm}\label{thm:pathologychiplusdelta}
    $P(\chi+\delta)=2$. 
\end{thm}
The proofs of our results follow in the next sections.

\section{Proof of Theorem \ref{thm:bigcapmatrix}}

Before we proceed to the proof of Theorem \ref{thm:bigcapmatrix}, we recall some definitions about trees. 
Denote by $\omega^{<\omega}$ the set of finite sequences with values in $\omega$. We say that $a=(a_0,\ldots,a_n) \in \omega^{<\omega}$ is an extension of $b=(b_0,\ldots,b_m) \in \omega^{<\omega}$, shortened as $b\subseteq a$, if $n\ge m$ and $a_k=b_k$ for all $k\le m$. Given an infinite sequence $z \in \omega^{\omega}$, we define $z\upharpoonright n:=(z_0,\ldots,z_{n-1}) \in \omega^{<\omega}$ for each $n \in \omega$. A \emph{tree on} $\omega$ is a subset $T\subseteq \omega^{<\omega}$ with the property that $a \in T$ for all $a,b \in \omega^{<\omega}$ such that $a\subseteq b \in T$. The \emph{body} of a tree $T$ on $\omega$, denoted by $[T]$, is the set of all its infinite branches, that is, the set of all sequences $z \in \omega^{\omega}$ such that $z\upharpoonright n \in T$ for all $n \in \omega$. 

Let $\mathrm{Tr}$ be the set of trees on $\omega$, see \cite[Definition 2.1]{K}. 
Identifying each tree with its characteristic function, the set $\mathrm{Tr}$ can be regarded as a closed subset of the Polish space $\{0,1\}^{\omega^{<\omega}}$, see \cite[Exercise 4.32]{K}. Thanks to Alexandrov's theorem \cite[Theorem 3.11]{K}, $\mathrm{Tr}$ is a Polish space. 
Now, let $\mathrm{WF}$ be the set of well-founded trees on $\omega$, that is, 
$$
\mathrm{WF}:=\left\{T \in \mathrm{Tr}: [T]=\emptyset\right\}.
$$
It is known that $\mathrm{WF}$ is a $\Pi^1_1$-complete subset of $\mathrm{Tr}$, hence it is coanalytic but not analytic, see \cite[Theorem 27.1]{K}. In particular, $\mathrm{WF}$ is not analytic in $\mathrm{Tr}$.

\begin{proof}
[Proof of Theorem \ref{thm:bigcapmatrix}] 
In this proof, we fix a bijection $f:\omega^{<\omega}\to\omega$. 
For each $z\in\omega^{\omega}$, define the infinite matrix $A^z=\left(a^z_{n,k}: n,k\in\omega\right)$ by setting 
$$
a^z_{n,k}
        :=
        \begin{cases}
        \,1 & \text{if } k=f(z\upharpoonright (n+1)),\\
        \,0 & \text{otherwise},
        \end{cases}
$$
for all $n,k \in \omega$. 
Observe that each matrix $A^z$ is nonnegative regular. Indeed, it is clear that $\sum_{k}a^z_{n,k}=1$ for all $n \in \omega$. In addition, if $k\in \omega$ is fixed, then $a^z_{n,k}=1$ for at most one $n\in\omega$, hence $\lim_{n}a^z_{n,k}=0$. To complete the proof, it is enough to show that 
$$
\mathcal I:=\bigcap_{z\in\omega^{\omega}}\mathcal I_{A^z}
$$
is not analytic. 

To this end, pick an arbitrary $S\subseteq \omega$, and observe that 
\begin{displaymath}
    \begin{split}
        S \in \mathcal{I} \quad 
        & \Longleftrightarrow \quad 
        \forall z \in \omega^\omega,\,\, S \in \mathcal{I}_{A^z} \\
         & \Longleftrightarrow \quad 
        \forall z \in \omega^\omega,\,\, \lim_{n\to \infty}\sum_{k\in S}a_{n,k}^z=0\\
        & \Longleftrightarrow \quad 
        \forall z \in \omega^\omega,\exists n_z \in \omega, \forall n \ge n_z, \forall k \in S, \,\, a_{n,k}^z=0\\
        & \Longleftrightarrow \quad 
        \forall z \in \omega^\omega,\,\,\{n\in\omega:f(z\upharpoonright n)\in S\}\in\mathrm{Fin}.
    \end{split}
\end{displaymath}

Next, define the map 
$\Phi:\mathcal P(\omega^{<\omega})\to\mathcal P(\omega)$ by 
$$
  \forall A\subseteq  \omega^{<\omega}, \qquad 
  \Phi(A):=f[A].
$$
It is routine that $\Phi$ is a homeomorphism and, by the above observation, that 
\begin{equation}\label{eq:kdsjhgksrhg}
  \Phi(A) \in \mathcal{I} 
  \,\,\,\Longleftrightarrow \,\,\, 
  \forall z\in\omega^{\omega},\,
        A\cap\{z\upharpoonright n:n\in\omega\}\in\mathrm{Fin}
\end{equation}
for all $A\subseteq  \omega^{<\omega}$. 

\begin{claim}\label{claim:kdsjhgsdg}
    $\Phi^{-1}[\mathcal I]\cap\mathrm{Tr}=\mathrm{WF}$. 
\end{claim}
\begin{proof}
Fix $T\in\mathrm{Tr}$. Suppose first that $T\notin\mathrm{WF}$.
Then it is possible to pick $z\in[T]$. Hence $z\upharpoonright n\in T$ for all $n\in\omega$, and therefore $T\cap\{z\upharpoonright n:n\in\omega\}$ is infinite. It follows by \eqref{eq:kdsjhgksrhg} that $\Phi(T)\notin\mathcal I$.

Conversely, suppose that $T\in\mathrm{Tr}$ and $\Phi(T)\notin\mathcal I$.
Again by \eqref{eq:kdsjhgksrhg}, there exists $z\in\omega^{\omega}$ such that $T\cap\{z\upharpoonright n:n\in\omega\}$ is infinite. Hence, for every $m\in\omega$, there exists $n\ge m$ such that
$z\upharpoonright n\in T$. Since $T$ is a tree, we obtain $z\upharpoonright m\in T$. By the arbitrariness of $m\in \omega$, it follows that $z\in [T]$. Thus
$T\notin\mathrm{WF}$.
\end{proof}

At this point, suppose for the sake of contradiction that $\mathcal I$ is analytic. 
Since $\Phi$ is a homeomorphism, $\Phi^{-1}[\mathcal I]$ is analytic in
$\mathcal P(\omega^{<\omega})$, see e.g. \cite[Proposition 14.4]{K}. Since $\mathrm{Tr}$ is closed, it follows
by Claim \ref{claim:kdsjhgsdg} that $\mathrm{WF}$ is analytic in $\mathrm{Tr}$, which is the desired contradiction. This concludes the proof. 
\end{proof}

\begin{rmk}
    It is routine to check that the ideal $\mathcal{I}$ defined in the proof of Theorem \ref{thm:bigcapmatrix} is not $P$-ideal. Since countable intersection of $P$-ideals is a $P$-ideal, it follows that $\mathcal{I}$ cannot be represented as a countable intersection of summable ideals. 
\end{rmk}

\section{Proof of Theorem \ref{thm:representation} and Corollary \ref{cor:nonpathological}} 

In order to prove Theorem \ref{thm:representation}, we need two intermediate lemmas. 
\begin{lem}\label{lem:fdsjgh}
    Let $A=(a_{n,k}: n,k \in \omega)$ be a nonnegative regular matrix and, for each $n \in \omega$ and $S\subseteq \omega$, set $\mu_n(S):=\sum_{k \in S}a_{n,k}$. Pick also a compact $\mathcal{K}\subseteq \mathcal{I}_A$, an infinite $L\subseteq \omega$, and a positive $\varepsilon>0$. Then, for each $n_0 \in \omega$, there exist a nonempty finite $F\subseteq L$ with $\min F\ge n_0$ and nonnegative reals $(c_n: n \in F)$ such that 
    \begin{equation}\label{eq:claimlemma}
    \sum_{n \in F}c_n=1 
    \quad \text{ and }\quad 
    \sup_{S \in \mathcal{K}}\sum_{n \in F}c_n \mu_n(S)<\varepsilon.
    \end{equation}
\end{lem}
\begin{proof}
Fix $n_0\in\omega$. 
Since $A$ is a nonnegative regular matrix, it follows that $M:=\sup_n \mu_n(\omega)<\infty$ and $\lim_n a_{n,k}=0$ for all $k \in \omega$. For each $n \in \omega$, pick $t(n) \in \omega$ such that $\sum_{k\geq t(n)}a_{n,k}<\varepsilon/2$. Thus, define a function \(f_n:\mathcal K\to[0,\infty)\) by
\[
\forall S\in \mathcal{K}, \qquad 
    f_n(S):=\sum_{k\in S, k<t(n)}a_{n,k}.
\]
Each $f_n$ is continuous, and $0\le f_n(S) \le M$ for all $n \in \omega$ and $S \in \mathcal{K}$. 
Since $\mathcal{K}\subseteq \mathcal{I}_A$, we get $\lim_n f_n(S)=0$ for all $S \in \mathcal{K}$.  Hence \((f_n)\) is a bounded sequence in the
Banach space \(C(\mathcal K)\), equipped with the supremum norm, and it
converges pointwise to \(0\).

\smallskip

We claim that \(f_n\to 0\) weakly in \(C(\mathcal K)\). Indeed, pick a bounded linear operator $T: C(\mathcal K)\to \mathbb{R}$. 
By the Riesz representation theorem, there is
a finite signed regular Borel measure \(\sigma\) on \(\mathcal K\) such that 
$T(g)=\int_{\mathcal K}g\,\mathrm{d}\sigma$ for all $g \in C(\mathcal{K})$, see e.g. \cite[Examples 1.2.11 and 1.10.6]{MR1650235}.  
Since \(|f_n|\leq M\) and
\(f_n(S)\to 0\) for every \(S\in\mathcal K\), the dominated convergence
theorem gives
\[
    \lim_{n\to \infty}T(f_n)=\lim_{n\to \infty}\int_{\mathcal K}f_n\,\mathrm{d}\sigma
    =\int_{\mathcal K}\lim_{n\to \infty}f_n\,\mathrm{d}\sigma=0.
\]
Thus \(f_n\to 0\) weakly in \(C(\mathcal K)\).

\smallskip

Now, let $(\ell_i: i \in \omega)$ be the increasing enumeration of the infinite set \(L\cap[n_0,\infty)\). 
Then \((f_{\ell_i})_{i\in\omega}\) is also weakly null in \(C(\mathcal K)\).
By Mazur's lemma, see e.g. \cite[Theorem 2.5.16]{MR1650235}, \(0\) belongs to the norm-closed convex hull of every tail
of a weakly null sequence. 
Therefore there are indices
\(i_0,\ldots,i_r\in\omega\) and nonnegative reals
\(d_0,\ldots,d_r\) such that
\[
    \sum_{q=0}^r d_q=1
    \qquad\text{and}\qquad
    \left\|\,\sum_{q=0}^r d_q f_{\ell_{i_q}}\,\right\|_\infty
    <\frac{\varepsilon}{4}.
\]

\smallskip

Thus, define $F:=\{\ell_{i_0},\ldots,\ell_{i_r}\}\subseteq L\cap[n_0,\infty)$. 
It follows that we may write the above convex combination as $\sum_{n\in F}c_n f_n$, where \(F\) is a nonempty finite subset of \(L\), \(\min F\geq n_0\),
\(c_n\geq 0\) for all $n \in F$, and \(\sum_{n\in F}c_n=1\). Moreover, the latter inequality can be rewritten as 
$$
\sup_{S\in\mathcal K}\sum_{n\in F}c_n f_n(S)<\frac{\varepsilon}{4}.
$$

Finally, for each \(S\in\mathcal K\), we have 
\[
\begin{aligned}
    \sum_{n\in F}c_n\mu_n(S)
    &=
    \sum_{n\in F}c_n
    \left(
        \sum_{\substack{k<t(n)\\ k\in S}}a_{n,k}
        +
        \sum_{\substack{k\geq t(n)\\ k\in S}}a_{n,k}
    \right)                                      \\
    &\leq
    \sum_{n\in F}c_n f_n(S)
    +
    \sum_{n\in F}c_n\sum_{k\geq t(n)}a_{n,k}      <
    \frac{\varepsilon}{4}
    +
    \sum_{n\in F}c_n\frac{\varepsilon}{2}=\frac{3\varepsilon}{4}.
\end{aligned}
\]
We conclude that \eqref{eq:claimlemma} holds. 
\end{proof}

\medskip

\begin{lem}\label{lem:nonpathlscsm}
    Let $\lambda: \mathcal{P}(\omega)\to [0,\infty]$ be a $\sigma$-additive measure such that $\lambda(\{k\})<\infty$ for all $k \in \omega$. Then 
    $
    \lambda \wedge k
    $ 
    is a nonpathological lscsm for each $k \in \omega$. 
\end{lem}
\begin{proof}
If $\lambda=0$ the claim is obvious. Hence, suppose hereafter that $\lambda \neq 0$. 
It follows by \cite[Proposition 4.1.6(a)--(b)]{FTpathology} that $\lambda$ is a nonpathological lscsm. 
Hence it is routine to show that $\lambda \wedge k$ is a nonpathological lscsm as well for each $k \in \omega$. 
\end{proof}

\medskip

Finally, we proceed to the proof of Theorem \ref{thm:representation}. 

\begin{proof}
    [Proof of Theorem \ref{thm:representation}] 
    \ref{item:1repres} $\implies$ \ref{item:2repres}. See \cite[Theorem 11.2]{FTpathology}. 

\medskip

    \ref{item:2repres} $\implies$ \ref{item:3repres}. As remarked in \cite[p. 23]{FTpathology}, this follows by the fact that every summable ideal can be represented as intersection of matrix summability ideals; cf. \cite[Theorem 8.5]{FTpathology}.

\medskip

    \ref{item:3repres} $\implies$ \ref{item:1repres}. Let $(A_j: j \in \omega)$ be a family a nonnegative regular matrices such that 
    $$
    \mathcal{I}:=\bigcap_{j \in J}\mathcal{I}_{A_j}
    $$
    is an $F_\sigma$ ideal. We need to show that there exists a nonpathological lscsm $\varphi: \mathcal{P}(\omega)\to [0,\infty]$ such that $\mathcal{I}=\mathrm{Fin}(\varphi)$. 

    For each $n\in \omega$ and $j \in J$, define the $\sigma$-additive measure $\mu^j_n: \mathcal{P}(\omega) \to [0,\infty]$ by 
    $$
\forall S\subseteq \omega, \qquad    \mu^j_n(S):=\sum_{k \in S}a^j_{n,k}. 
    $$
    Observe that, for each $j \in J$, the regularity of $A_j$ implies $\sup_n \mu^j_n(\omega)<\infty$. 

    Since $\mathcal{I}$ is an $F_\sigma$ ideal, it is possible to choose an increasing sequence $(\mathcal{K}_m: m \in \omega)$ of subsets of $\omega$ such that each $\mathcal{K}_m$ is hereditary compact and, in addition, 
    $
    \{\{n\}: n \in \omega\}\subseteq K_0.
    $ 
    Now, fix $S\notin\mathcal I$ and $p\in \omega$. 
    Pick an index $j \in J$ (depending on $S$) such that $S\notin\mathcal I_{A_j}$, so that $\limsup_n\mu^j_n(S)>0$. Choose $\varepsilon>0$ and an infinite set
$L\subseteq\omega$ such that 
$$
\forall n \in L, \qquad 
\mu^j_n(S)>\varepsilon.
$$
For every $m\geq p$, apply Lemma \ref{lem:fdsjgh} to the compact set $\mathcal{K}_m$,
the matrix $A_j$, the infinite set $L$, and $\eta=2^{-(m+1)}$. We obtain a finite set $F_m\subseteq L$ and coefficients
$(c^m_n)_{n\in F_m}$ such that $c^m_n\geq 0$, $\sum_{n\in F_m}c^m_n=1$, and
the measure
\[
    \nu_m:=\sum_{n\in F_m}c^m_n\mu^j_n
\]
satisfies $\nu_m(Y)\leq 2^{-(m+1)}$ for every $Y\in \mathcal{K}_m$. Since every $n\in F_m$ belongs to $L$, we also have $\nu_m(S)>\varepsilon$. At this point, define the map $\lambda_{S,p}: \mathcal{P}(\omega)\to [0,\infty]$ by 
$$
\forall Y\subseteq \omega, \qquad 
\lambda_{S,p}(Y):=\sum_{m\geq p}\nu_m(Y).
$$
It follows by construction that, if $Y \in \mathcal{K}_p\subseteq \mathcal{K}_m$ (hence, in particular, if $Y=\{k\}$ for some $k \in \omega$), then $\lambda_{S,p}(Y)\le \sum_{m\geq p}2^{-m-1} \le 1$. In addition, $\lambda_{S,p}$ is a $\sigma$-additive measure and $\lambda_{S,p}(Y)\ge \sum_{m\ge p}\varepsilon=\infty$. 

\smallskip

At this point, define the submeasure 
$$
\varphi:=\sup\left\{\lambda_{S,p}\wedge p: S\notin \mathcal{I}, p\in \omega\right\}.
$$
As it follows by Lemma \ref{lem:nonpathlscsm}, each $\lambda_{S,p}\wedge p$ is a nonpathological lscsm. As observed in \cite[Section 2.2]{FTpathology} $\varphi$ is a lscsm. In addition, $\varphi$ is nonpathological: indeed, pick $Y\subseteq \omega$ and let $(S_\alpha,p_\alpha)$ be a net for which $\varphi(Y)=\sup_\alpha (\lambda_{S_\alpha,p_\alpha}\wedge p_\alpha)(Y)$. Taking into account that each $\lambda_{S,p}\wedge p$ is nonpathological, it follows that 
$$
\varphi(Y)=\sup_\alpha \sup\{\mu(Y): \mu \text{ is a finitely additive measure with }\mu\le \lambda_{S_\alpha,p_\alpha}\wedge p_\alpha\}.
$$
This implies that $\varphi(Y)=\widehat{\varphi}(Y)$. Therefore $\varphi$ is a nonpathological lscsm such that $\varphi(\{k\})\le 1$ for all $k \in \omega$. 

\smallskip

To complete the proof, let us show finally that 
$$
\mathcal{I}=\mathrm{Fin}(\varphi).
$$
First, pick $Y\in\mathcal I$ and choose an integer $m\ge 1$ such that $Y\in \mathcal{K}_m$. If $p\geq m$, then
$Y\in \mathcal{K}_m\subseteq \mathcal{K}_p$, so
$\lambda_{S,p}(Y)\wedge p\le \lambda_{S,p}(Y)\leq 1$ for every $S\notin\mathcal I$. If $p<m$, then trivially
$\lambda_{S,p}(Y)\wedge p\leq p<m$ for every $S\notin \mathcal{I}$. It follows that $\varphi(Y)\leq m<\infty$. This proves that $\mathcal I\subseteq\operatorname{Fin}(\varphi)$. 

Conversely, pick a set $Y\notin\mathcal I$. Since  $\lambda_{Y,p}(Y)=\infty$, we obtain that 
$$
\varphi(Y)\ge \sup\{(\lambda_{Y,p}\wedge p)(Y): p \in \omega\}=\sup\{p: p \in \omega\}=\infty.
$$
Hence $Y\notin \mathrm{Fin}(\mathcal{I})$. Therefore $\operatorname{Fin}(\varphi)\subseteq \mathcal{I}$. This completes the proof. 
\end{proof}

\medskip

\begin{proof}
    [Proof of Corollary \ref{cor:nonpathological}]
    Thanks to \cite[Corollary 3.12]{MR4756411}, $\mathcal{I}$ is a nonpathological analytic $P$-ideal $\mathcal{I}$ on $\omega$ if and only if it can be represented as intersection of matrix summability ideals. At the same time, it follows by Theorem \ref{thm:representation} that the latter condition is equivalence to the fact that $\mathcal{I}$ is a nonpathological $F_\sigma$ ideal. 
\end{proof}

\section{Proof of Theorem \ref{thm:randomgraph}}

Let us recall that, if $(\omega,E)$ is the random graph, then for all disjoint finite subsets $F,G\subseteq \omega$, there exists $k \in \omega\setminus (F\cup G)$ such that $\{\{n,k\}: n \in F\}\subseteq E$ and $\{\{n,k\}: k \in G\}\cap E=\emptyset$. In addition, given a countable graph $G$ with edges $D\subseteq [G]^2$, there exists $S\subseteq \omega$ such that $(G,D)$ is isomorphic to $(S,E\upharpoonright S)$, see e.g. \cite[Lemma 1.6.28 and Lemma 1.6.29]{Mezathesis}. 
\begin{proof}
    [Proof of Theorem \ref{thm:randomgraph}] 
    First, it is known that $\mathcal{R}$ is an $F_\sigma$ ideal, see e.g. \cite[p.30]{Mezathesis}. To complete the proof, suppose for the sake of contradiction that there exists a nonpathological lscsm $\varphi$ such that $\mathcal{R}=\mathrm{Fin}(\varphi)$. 

    For each integer $m\geq 1$, let $\mathsf{KG}(3m,m)$ be the Kneser graph whose vertices are the $m$-element subsets of $\{1,2,\ldots,3m\}$, where two vertices are adjacent if and only if the two corresponding sets are disjoint. 
By the Kneser--Lovasz theorem \cite{MR514625}, the chromatic number of $\mathsf{KG}(3m,m)$ (that is, the smallest number of colors needed for a labeling of the vertices such that no two vertices sharing the same edge have the same color)
is exactly 
$$
\mathrm{Chrom}(\mathsf{KG}(3m,m))=3m-2m+2=m+2.
$$
Since the random graph $(\omega,E)$ contains every countable graph, we may pick pairwise disjoint finite sets $(B_m: m\geq 1)$ such that $\omega\upharpoonright B_m$ is isomorphic to $\mathsf{KG}(3m,m)$ for each $m\geq 1$, and there are no edges between $B_m$ and $B_n$ for all distinct $m,n \in \omega$. 

\begin{claim}\label{claim:R1}
    For each $N\ge 1$, we have $T_N:=\bigcup_{m\ge N}B_m\notin \mathcal{R}$. 
\end{claim}
\begin{proof}
    Fix $N\ge 1$ and suppose for the sake of contradiction that $T_N$ is covered by finitely many cliques and anticliques of the random graph, let us say by $q$ cliques and $r$ anticliques. Since there are no edges between
distinct blocks $B_m$, each of the $q$ cliques meets at most one block. Hence, for all but finitely many $m\geq N$, the whole block $B_m$ is covered by the $r$ anticliques. This implies that $\mathrm{Chrom}(\omega \upharpoonright B_m)\leq r$ for all
but finitely many $m\geq N$, contradicting
$\mathrm{Chrom}(\omega\upharpoonright B_m)=m+2$.
\end{proof}

Let $\mathcal{K}$ be the family of all anticliques contained in $T_1$. 
We regard $\mathcal{K}$ as a subspace of $\mathcal{P}(T_1)$. 
Then $\mathcal{K}$ is a compact subset of $\mathcal{R}=\mathrm{Fin}(\varphi)$. 
Hence $\mathcal{K}=\bigcup_j \mathcal{K}_j$, where 
\[
\forall j \in \omega, \qquad         \mathcal{K}_j:=
        \{A\in\mathcal{K}:\varphi(A)\leq j\}.
\]
Observe that each $\mathcal{K}_j$ is closed in $\mathcal{K}$ by lower
semicontinuity of $\varphi$. Hence, by the Baire category theorem, 
there exist
an integer $M\in\omega$ and finite sets $s,t\subseteq T_1$ such that
$\varphi(A)\leq M$ for every $A\in\mathcal{K}$ satisfying $s\subseteq A$ and $A\cap t=\emptyset$.
Choose $N\geq 1$ such that $s\cup t\subseteq \bigcup_{m<N}B_m$. Since there are
no edges between distinct blocks, it follows that $s\cup J\in\mathcal{K}$ for
every anticlique $J\subseteq T_N$. In addition, $(s\cup J)\cap t=\emptyset$.
Therefore $\varphi(J)\leq \varphi(s\cup J)\leq M$ for every anticlique $J\subseteq T_N$.

\begin{claim}\label{claim:R2}
    $\varphi(F)\leq 3M$ for all finite $F\subseteq T_N$. 
\end{claim}
\begin{proof}
Fix a finite set $F\subseteq T_N$, and let
$\mu:\mathcal{P}(\omega)\to[0,\infty]$ be a finitely additive measure such that
$\mu\leq\varphi$. Let
\[
        Q:=\{m\geq N:F\cap B_m\neq\emptyset\}.
\]
Then $Q$ is finite. Fix $m\in Q$. For each $i\in\{1,\ldots,3m\}$, define
\[
        S_i^m:=\{a\in B_m:i\in a\}.
\]
Each $S_i^m$ is an anticlique of $\omega\upharpoonright B_m$, since two $m$-element subsets of $\{1,\ldots,3m\}$ containing $i$ cannot be disjoint.
Moreover, every vertex of $B_m$ belongs to exactly $m$ of the sets
$S_1^m,\ldots,S_{3m}^m$. Since $F\cap B_m$ is finite, finite additivity gives
\[
        \sum_{i=1}^{3m}\mu(F\cap S_i^m)
        =
        m\,\mu(F\cap B_m).
\]
Hence there exists $i_m\in\{1,\ldots,3m\}$ such that
\[
        \mu(F\cap S_{i_m}^m)
        \geq
        \frac{1}{3}\mu(F\cap B_m).
\]
Set
\[
        J:=\bigcup_{m\in Q}(F\cap S_{i_m}^m).
\]
Then $J$ is an anticlique contained in $T_N$, so that $\varphi(J)\le M$ by the above observation. Since $\mu\leq\varphi$, we
obtain
\[
\begin{split}
        \mu(F)
        &=
        \sum_{m\in Q}\mu(F\cap B_m)                                      \leq
        3\sum_{m\in Q}\mu(F\cap S_{i_m}^m)                                =
        3\mu(J)
        \leq
        3\varphi(J)
        \leq 3M.
\end{split}
\]
Since this holds for every finitely additive measure $\mu\leq\varphi$, and
$\varphi$ is nonpathological, we conclude that $\varphi(F)=\widehat{\varphi}(F)\le 3M$. 
\end{proof}

It follows by Claim \ref{claim:R2} and the lower semicontinuity of $\varphi$ that 
$
        \varphi(T_N)
        =
        \lim_{n}\varphi(T_N\cap n)
        \leq 3M
        <\infty.
$
Thus $T_N\in\mathrm{Fin}(\varphi)=\mathcal{R}$. This contradicts Claim \ref{claim:R1} and provides the desired contradiction, concluding the proof.
\end{proof}

\section{Proof of Theorem \ref{thm:pathologychiplusdelta}} 

Before we proceed to the proof of Theorem \ref{thm:pathologychiplusdelta}, we recall some basic properties about the submeasure $\chi$ defined in \eqref{eq:definitionchi}.
\begin{lem}\label{lem:basicchi}
    The following properties hold\textup{:}
    \begin{enumerate}[label={\rm (\roman*)}]
    \item \label{item:1basic} $\chi_n(A(n,I))\leq 2^{1-n}|I|$ for each $n\ge 1$ and $I\subseteq \Delta_n$\textup{;}
    \item \label{item:2basic} $\chi_n(X_n)=1+2^{1-n}$ for all $n\ge 1$, hence $\chi(X)=1$\textup{;}
    \item \label{item:3basic} $\chi(F)=0$ for all finite $F\subseteq X$\textup{.}
    \end{enumerate}
\end{lem}
\begin{proof}
     Item \ref{item:1basic} follows by definition; the proof of item \ref{item:2basic} can be found in \cite{MR419712}. Item \ref{item:3basic} follows by the fact that the superior limit in \eqref{eq:definitionchi} ignores finite sets. 
\end{proof}

Also, we will use the following key property, which will be proved with the probabilistic method. 
\begin{prop}\label{prop:key}
    For each $\varepsilon \in (0,1)$ and for all integers $r,L\ge 1$ such that $r\le \varepsilon L$ and $(1-\varepsilon/2)^{(r+1)/2}<\nicefrac{1}{2}$, there exists an integer $n_0=n_0(\varepsilon,r,L) \ge 1$ such that for all $n\ge n_0$ there are sets 
    $$
    I_{n,1},\ldots,I_{n,L}\subseteq \Delta_n
    $$
    for which $|I_{n,j}|\le \varepsilon 2^n$ for all $j\le L$ and $|\{j\le L: x \cap I_{n,j}=\emptyset\}|\le r$ for all $x\in X_n$. 
\end{prop}
\begin{proof}
Choose $I_{n,1},\ldots,I_{n,L}$ randomly and independently by putting each element of $\Delta_n$ into each $I_{n,j}$ independently with probability $p:=\varepsilon/2$. Define
$$
q:=\left(1-p\right)^{(r+1)/2} \in \left(0,\,\frac{1}{2}\right).
$$

Fix $x\in X_n$ and $J\subseteq\{1,\dots,L\}$ with $|J|=r+1$. Taking into account that $|x|=2^{n-1}$, the probability that $x$ misses every $I_{n,j}$ for $j\in J$ is
$$
\left(1-p\right)^{(r+1)2^{n-1}}=q^{2^n}.
$$

Since $|X_n|\le |\mathcal{P}(\Delta_n)|=2^{2^n}$, it follows that the probability that there exist
$x\in X_n$ and $J\subseteq\{1,\dots,L\}$ with $|J|=r+1$ such that $x\cap I_{n,j}=\emptyset$ for every $j\in J$ is at most
$$
\binom{L}{r+1}
2^{2^n}
\left(1-p\right)^{(r+1)2^{n-1}}
\le 2^L \cdot 
(2q)^{2^n}.
$$
This tends to $0$ as $n\to\infty$, because $2q \in (0,1)$ by hypothesis. 

On the other hand, for each fixed $j$, the random variable $|I_{n,j}|$ has
binomial distribution with mean $2^np$. 
Hence by the weak law of large numbers 
$$
\lim_{n\to \infty}\mathbb P\bigl(|I_{n,j}|>\varepsilon 2^n\bigr)=0. 
$$
Since $L$ is fixed, with probability tending to \(1\) we
have $|I_{n,j}|\leq \varepsilon 2^n$ for every \(j=1,\dots,L\).

Thus, for all sufficiently large \(n\), there exists a choice of
\(I_{n,1},\dots,I_{n,L}\) satisfying both required properties.
\end{proof}

We conclude with the proof of Theorem \ref{thm:pathologychiplusdelta}. Here, for each $x \in X$, let $e_x: \mathcal{P}(X)\to [0,\infty]$ be the submeasure defined by 
$$
\forall B\subseteq X, \qquad 
e_x(B)=
\begin{cases}
\,0 & \text{ if } x\notin B,\\
\,1 & \text{ if }x \in B.
\end{cases}
$$
\begin{proof}
     [Proof of Theorem \ref{thm:pathologychiplusdelta}]
     Set $\varphi:=\chi+\delta$. First, we start proving that 
     \begin{equation}\label{eq:importantclaim}
     \widehat{\varphi}(X)=1.
     \end{equation}
     In fact, for each $x \in X$ we have that $e_x$ is a finitely additive measure such that $e_x\le \varphi$. This implies that $\widehat{\varphi}(X)\ge e_x(X)=1$. Viceversa, let $\mu: \mathcal{P}(X)\to [0,\infty]$ be a finitely additive measure such that $\mu\le \varphi$. It will be enough to show that $\mu(X)\le 1$. To this aim, fix $\varepsilon \in (0,1)$ and pick two integers $r,L\ge 1$ such that $r\le \varepsilon L$ and $(1-\varepsilon/2)^{(r+1)/2}<\nicefrac{1}{2}$. Hence, pick $n_0 \ge 1$ and, for each $n\ge n_0$, choose sets $I_{n,1},\ldots,I_{n,L}$ satisfying the claim of Proposition \ref{prop:key}. 
     For each $j=1,\dots,L$, define the set 
$$
B_j
=
\left(\bigcup_{n<n_0}X_n\right)
\cup
\left(\bigcup_{n\geq n_0} A(n,I_{n,j})\right).
$$
Since the finite initial part $\bigcup_{n<n_0}X_n$ does not affect the superior limit in \eqref{eq:definitionchi}, we have
$$
\chi(B_j)
=
\limsup_{n\to\infty}\chi_n(A(n,I_{n,j})).
$$
Thanks to Lemma \ref{lem:basicchi}\ref{item:1basic} and Proposition \ref{prop:key}, we have 
$$
\chi_n(A(n,I_{n,j}))\le 2^{1-n}|I_{n,j}|\le 2^{1-n}(\varepsilon 2^n)=2\varepsilon
$$
for all $n\ge n_0$ and all $j=1,\ldots,L$. This implies that 
$$
\mu(B_j)\le \varphi(B_j)=\chi(B_j)+\delta(B_j)\le 1+2\varepsilon
$$
for all $j=1,\ldots,L$.

Now we use the covering property of the $B_j$'s. Every point of the finite
initial part 
$\bigcup_{n<n_0}X_n$ 
belongs to all $B_j$'s. If $x\in X_n$ for $n\geq n_0$, it follows by Proposition \ref{prop:key} that 
$x$ misses at most $r$ of the sets $I_{n,j}$. 
Hence $x$ belongs to at least $L-r$ (which is at least $L-\varepsilon L$) of the sets $B_1,\dots,B_L$. 
Therefore, pointwise on $X$, we obtain
$$
\mathbf 1_X
\leq
\sum_{j=1}^L \frac{1}{(1-\varepsilon)L}\mathbf 1_{B_j}.
$$

At this point, for each $S\subseteq\{1,\dots,L\}$, define 
$$
C_S=
\{x\in X: x\in B_j\ \text{if and only if}\ j\in S\}.
$$
Then the sets $C_S$, some of which may be empty, form a finite partition of $X$. 
Moreover, if $C_S\neq\emptyset$, then $|S|\geq (1-\varepsilon)L$. 
Hence, using finite additivity of $\mu$, we get 
\begin{displaymath}
    \begin{split}
\mu(X)
&=
\sum_{S\subseteq\{1,\dots,L\}}\mu(C_S)\leq
\sum_{S\subseteq\{1,\dots,L\}}
\frac{|S|}{(1-\varepsilon)L}\mu(C_S)\\
&=\sum_{j=1}^L \frac{1}{(1-\varepsilon)L}\mu(B_j)
\le \sum_{j=1}^L \frac{1}{(1-\varepsilon)L}(1+2\varepsilon)
=\frac{1+2\varepsilon}{1-\varepsilon}.
    \end{split}
\end{displaymath}
Since $\varepsilon>0$ was arbitrary, we conclude that $\mu(X)\le 1$. 

\medskip

Now, pick a nonempty $B\subseteq X$, a point $x \in B$, and a finitely additive measure $\mu\le \varphi$. It follows by \eqref{eq:importantclaim} that
$$
1=e_x(B)\le \widehat{\varphi}(B) \le \widehat{\varphi}(X)=1,
$$
This implies that 
$$
\widehat{\varphi}=\delta.
$$
Therefore 
$$
P(\varphi)=\sup\left\{\frac{\varphi(B)}{\delta(B)}: B\subseteq X\right\}=
1+\sup\left\{\frac{\chi(B)}{\delta(B)}: B\neq \emptyset\right\}=1+\chi(X)=2,
$$
where the last equality follows by Lemma \ref{lem:basicchi}\ref{item:2basic}. 
\end{proof}

\bibliographystyle{amsplain}

\end{document}